\UseRawInputEncoding
\documentclass[10pt, a4paper]{article}
\usepackage[latin1]{inputenc}

\usepackage{amsmath,amssymb,amsthm,color,physics}
\usepackage{enumitem}
\newtheorem*{assumption*}{\assumptionnumber}
\providecommand{\assumptionnumber}{}
\makeatletter
\newenvironment{assumption}[1]
 {%
  \renewcommand{\assumptionnumber}{Hypothesis $\left(#1\right)$}%
  \begin{assumption*}%
  \protected@edef\@currentlabel{$\left(#1\right)$}%
 }
 {%
  \end{assumption*}
 }
\usepackage{graphicx}
\usepackage[hang]{caption}

\def\X{{\mathcal{X}}}

\def\Vo{{\mathbb{V}}}

\def\ball{{I\kern -.35em B}}

\def\X{{\mathcal{X}}}

\def\ball{{I\kern -.35em B}}

\def\tto{\;{\lower 1pt \hbox{$\rightarrow$}}\kern -12pt
            \hbox{\raise 2.8pt \hbox{$\rightarrow$}}\;}
\def\dom{\mathop{\rm dom}\nolimits}

\def\bz{\bar z}

\def\bx{\bar x} 
\def\by{\bar y}

\def\R{{\mathbb R}}
\def\N{{\mathbb N}}

\def\dis{\displaystyle}
\def\inte{{\rm int}}
\def\cl{{\rm cl}}

\def\dom{\mathop{\rm dom}\nolimits}

\def\lev{\mathop{\rm lev}\nolimits}

\def\T{{\mathcal T}}

\def\F{{\mathcal{F}}}
\def\R{{\mathbb{R}}}

\def\P{{\mathcal{P}}}

\def\xnk{x_{n_k}}

\def\Se{{\mathcal{S}}{\rm Eff}}
\def\Re{{\mathcal{R}}{\rm Eff}}
\def\Pe{{\mathcal{P}}{\rm Eff}}
\def\Ge{{\mathcal{G}}{\rm Eff}}

\def\Te{{\mathcal{T}}{\rm Eff}}

\def\T{{\mathcal{T}}}

\def\topo{\tau^{C}}

\def\pzs{{\mathcal{P}}^\circ(Z)} 
\def\pxs{{\mathcal{P}}^\circ(X)}

\def\ppz{\P_p(Z)}
\def\ppR{\P_p(\R)}
\def\ppR2{\P_p(\R^2)}

\def\ll{{\,\prec_C^{l}\,}}

\def\lel {{\,\dis \preceq_{C}^l\,}}
\def\lell{{\,\dis \curlyeqprec_{C}^l\,}}
\def\nlell{{\,\dis \not \curlyeqprec_{C}^l\,}}

\def\Ns{{\mathcal{N}_\infty}}
\def\Nss{{\mathcal{N}^\#_\infty}}

\def\X{\mathcal{X}}

\def\CM{K_M}

\def\X{{\mathrm{X}}}
\def\ball{{I\kern -.35em B}}

\def\tto{\;{\lower 1pt \hbox{$\rightarrow$}}\kern -12pt
\hbox{\raise 2.8pt \hbox{$\rightarrow$}}\;}

\def\hX{\widehat{X}}

\def\hx{\hat{x}}

 \def\bz{\bar z} 
   \def\bx{\bar x} 
  \def\by{\bar y}

\def\R{{\mathbb R}}
\def\N{{\mathbb N}}

\def\F{\mathcal{G}_C\left(X, \ppz \right)}

\def\FF{\mathcal{G}_C\left(X, \pzs \right)}

\def\DNFF{\mathcal{G}_C\left(D_n, \pzs \right)}

\usepackage{cite}

\newtheorem{definition}{Definition}[section]
\newtheorem{theorem}[definition]{Theorem}
\newtheorem{lemma}[definition]{Lemma}
\newtheorem{remark}[definition]{Remark}
\newtheorem{proposition}[definition]{Proposition}

\newtheorem{example}[definition]{Example}

\title{Stability analysis for set-valued optimization in Geoffroy spaces}
 
\author{James Larrouy\footnote{LAMIA, Universit\'e des  Antilles, Pointe-\`a-Pitre, Guadeloupe (France), {\tt james.larrouy@univ-antilles.fr}}}

\date{$\star$}
\begin{document}
\maketitle

\noindent{\bf Abstract.} In this work, we study the external and internal stability of minimal solutions to set-valued optimization problems in a new functional framework. We consider perturbations on both the objective function and the admissible domain. To address these problems, we introduce two variational convergences for sequences of set-valued maps, namely the Gamma-cone convergence and the sequential Gamma-cone convergence. The upper and the lower convergence of strong level sets are also studied. \medbreak



\def\X{\mathcal{X}}

\def\topoK{\tau^{\CM}}
\def\bt{\overline{\theta}}
\def\nll{\nprec^l_{C}} 
\def\nlel{\npreceq^l_{C}} 


\section{Introduction}\label{sec:1}

Over the past decade, novel challenges have emerged within the domain of financial risk management, as well as in the modelling and treatment of high uncertainties (see for instance \cite{HHE10,HHR11,KTZ15}). The corresponding contemporary challenges generally lead to the study of optimization problems involving an objective function with values in a preordered hyperspace - by a convex cone. These are set-valued optimization problems where the objective functions usually have as codomain the set of the nonempty subsets of a real normed vector space. Based on the theory of conlinear spaces (see Hamel \cite{HAM05}), Michel H. Geoffroy \cite{GEO19} proposed in 2019 a topological convergence for sequences of sets adapted to the treatment of set optimization problems over the powerset constisting of lower bounded sets. More recently, Geoffroy and Larrouy \cite{GLA22} extended this work by proposing a topological framework for conlinear spaces, devoted to handling set-valued optimization problems, called set order topology. From now on, we will refer to any space of set-valued maps whose codomain is a conlinear space equipped with the set order topology as a set order functional space. 

In the context of set-valued optimization, stability analysis refers to the study of the asymptotic behavior of solution sets under data perturbation. In their pioneering work, Guti\'errez \textit{et al.} \cite{GMM16} introduced the concepts of external (or upper) and internal (or lower) stability for set-valued optimization problems based on the set approach \cite{KUR97, KUR98, KUR01}. Since then, many authors have taken an interest in this type of problem (see \cite{ADH20, ADH20b, GMN17, GEO19, GLA22, HAN22, HLO22, HLO22b, KLA19, KLA24, KHO18} and references therein), in particular because the set approach is meaningful not only in terms of mathematics but also in practical applications. In the case of external stability, the objective is to establish that any sequence of suitably chosen solution sets of perturbed problems converges to the solution set of the original problem, whereas in the case of internal stability, each solution set of the original problem is equivalently estimated by a subsequence of solution sets of perturbed problems.  It should be noted that virtually all recent work dealing with the stability of minimal solutions to set-valued optimization problems pays particular attention to the structure of the codomain of the objective function (often called image space). Furthermore, there is, in principle, no universal notion of stability in the literature. Most of the work on the topic has studied Painlev\'e-Kuratowski and/or Hausdorff convergence of sequences of solutions of approximated set optimization problems to solutions of the original problem in the image space (see \cite{ADH20, ADH20b, AHS22, GMM16, HAN22, KLA22, KLA24}). In their work, Karuna and Lalitha \cite{KLA19} introduced a notion of Gamma convergence to establish the external and internal stability in terms of Painlev\'e-Kuratowski convergence of sequence of solution sets of perturbed problems under certain compactness assumptions and domination properties. It is worth noting that in the context of vector optimization, several concepts of Gamma convergence have been proposed to address stability analysis (see for instance Oppezzi and Rossi \cite{ORO08, ORO08b} and Geoffroy \textit{et al.} \cite{GGJ17, GMN17}). \\

\noindent Therefore, since issues related to the structure of the image space have gained much attention in the litterature, we propose to study stability analysis in the context of set order topology, which is well suited when studying set-valued optimization problems. In this work, we introduce a new set order functional space called Geoffroy space along with the notion of Geoffroy minimality for set-valued optimization. Then, we propose a new variational convergence, namely the Gamma-cone convergence, which is helpful when it comes to investigate the asymptotic properties of level sets. We show that, in Geoffroy spaces, this concept is consistent with the well-known concept of Gamma convergence in the scalar case (see Braides \cite{BRA06}). Next, we extend in a natural fashion our concept of Gamma-cone convergence to the concept of sequential Gamma-cone convergence for sequences of set-valued maps whose domain depends on the index of each element. Then, we prove the external and the internal stability of strong, Pareto, Geoffroy and relaxed minimal solutions to set-valued optimization problems. Finally, we show how to use the results in practice and point out three open questions. \\

\noindent The paper is organized as follows: in Section \ref{Section 1}, we review some key concepts and provide some auxiliary results that will be used throughout. Then, we introduce the concept of Geoffroy spaces and the notion of Geoffroy minimality. In Section \ref{Section 2}, we first introduce the Gamma-cone convergence and we show that it possesses some essential properties. Afterward, we study the asymptotic properties of strong level sets from which we deduce the external and the internal stability of the sets of minimal solutions with respect to Painlev\'e-Kuratowski convergence. In Section \ref{Section 3}, we study the stability of minimal solutions in the sense of Guti\'errez \textit{et al.} \cite{GMM16}. We introduce the concept of Kuratoski pairs and of sequential Gamma-cone convergence which are used to study the external and the internal stability of minimal solution to set-valued optimization problems. Then, we show the practicability of the obtained results through a simple example.

\section{Preliminaries and notations}\label{Section Background materials}\label{Section 1}

In this section, we review some important notions and clarify the notation used in the paper. In addition, we provide some auxiliary results for the sequel and introduce a new functional structure that we will use in this work to deal with set-valued optimisation problems.

\subsection{Background materials}

Throughout this work, $Z := \left(Z, \Vert \cdot \Vert\right)$ stands for a nontrivial real normed vector space partially ordered by {a} nonempty closed, convex and solid cone $C \subsetneq Z$ through the binary relation $\leq_C$ defined by $z_1 \leq_C z_2$ if and only if $z_2 - z_1 \in C$. By cone we mean a set closed under nonnegative scalar multiplication. The notation $\inte(C)$ will be used to designate the topological interior of $C$ (with respect to the norm topology on $Z$). In this way, we can define the strict order binary relation $<_C$ on $Z$ by $y_1 \; <_C \; y_2$ if and only if $y_2 - y_1 \in \inte(C)$. On the other hand, the notation $\mathrm{cl}(\cdot)$ will be used to designate the topological closure of a given subset of $Z$. From now on, we denote $\pzs$ the set of all nonempty subsets of $Z$. We recall that the 
{$C$-lower} preorder $\lel$ {(or Kuroiwa's lower preorder)} is defined for $A, B \in \pzs$ by 
\begin{equation*}
A \; \lel B \; \Longleftrightarrow \;  B  \subseteq  A + C,
\end{equation*}
and that the 
{strict $C$-lower} preorder $\ll$ (see Geoffroy \textit{et al.} \cite{GLA22}) is defined by 
\begin{equation*}\label{StrictPreorder}
A \; \ll B \; \Longleftrightarrow  \; \exists {\varepsilon} \in \inte(C),\; A + {\varepsilon} \; \lel B.
\end{equation*}
Now, we introduce the large $C$-lower preorder defined by
\begin{equation*}\label{LargePreorder}
A \lell B \; \Longleftrightarrow  \; B  \subseteq  {\rm cl}(A + C),
\end{equation*}
and the associated large $C$-lower equivalence relation $\stackrel{l}{\approx}$ on $\pzs$ defined by 
\begin{equation*}\label{EquivalenceRelation}
A \stackrel{l}{\approx} B  \Longleftrightarrow A \; \lell B \text{ and  } B \; \lell A.
\end{equation*}
\noindent Obviously, $\ll$ forces $\lel$ and $\lel$ forces $\lell$. Moreover, these binary relations are compatible with the Minkowski sum and multiplication by positive scalars.

\begin{remark}\label{RemarkLinkRelations} Since $C$ is convex, $A \;  \ll B$ and $B \; \lell D$ forces $A \ll D$, for every $A, B, D \in \pzs$.
\end{remark}

\noindent  Before proceeding, we recall an important result of Tanaka and Kuroiwa \cite[Theorem 2.2]{TKU93}.

\begin{theorem}[\cite{TKU93}]\label{theoremTanakaKuroiwa}
Let $A, B \in \pzs$ be two convex sets. If $A$ is solid, then $\inte(A) + B = \inte(A + B).$
\end{theorem}

\noindent It follows from above that $\inte(C) + C = \inte(C).$ From now on, $\ball_{r}(z)$ stands for the open ball with center $z \in Z$ and radius $r > 0$. The result is as follows.

\begin{proposition}\label{PropositionInclusion} Let $A \in \pzs$ and $u \in \inte(C)$. Then, $\left\lbrace u \right\rbrace + \cl(A + C) \subset A + \inte(C).$
\end{proposition}

\begin{proof} Since $u \in \inte(C),$ there exists $r > 0$ such that $\dis \ball_r \left(\frac{u}{2}\right) \subset C$ which yields that $\dis \ball_r(u) \subset \inte(C).$ Now, take $\bz \in \cl(A + C).$ Then, from the first countability of $Z$, we can find a sequence $\left\lbrace z_n \right\rbrace \subset A + C$ converging to $\bz,$ and we can fix $n_0 \in \N$ such that $\Vert \bz - z_{n_{0}} \Vert < r.$ Let us consider the vector $z_0 =  \bz - z_{n_{0}} \in Z.$ Therefore, one can observe that $\Vert u + z_0 - u \Vert = \Vert z_0 \Vert < r.$ In other words, $u + z_0 \in \inte(C).$ Since $z_n \in A + C$ for all $n \in \N,$ we deduce from Theorem \ref{theoremTanakaKuroiwa} that $u + \bz = u + z_0 + z_{n_{0}} \in A + \inte(C).$
\end{proof}

\noindent It is worth noting that the strict and the large $C$-lower preorders are interrelated as follows.

\begin{lemma}\label{lemmastrictStrong} Let $A,B \in \pzs$. The following assertions hold true.
\begin{itemize}
\item[$(a)$] If $A + \varepsilon_0 \; \lell B$ for some $\varepsilon_0 \in \inte(C)$, then $A \;\ll B.$
\item[$(b)$] If $A-\varepsilon \; \ll B$ for all $\varepsilon \in\inte(C)$, then $A \;\lell B.$ 
\end{itemize}
\end{lemma}

\begin{proof} We prove $(a).$ Let us observe that since $Z$ is Hausdorff, $\cl (A + \varepsilon_0 + C) = \left\lbrace \varepsilon_0 \right\rbrace + \cl(A + C).$ Hence, using Proposition \ref{PropositionInclusion}, we get
\begin{equation*}
A+ \varepsilon_0 \; \lell B \; \Longleftrightarrow B \subset \left\lbrace \varepsilon_0 \right\rbrace + \cl(A + C) \; \Longrightarrow A + \frac{\varepsilon_0}{2} \; \lel B \; \Longrightarrow A \; \ll B.
\end{equation*}
Now, let us prove $(b)$. Take $\varepsilon \in \inte(C)$ and set $\varepsilon^n := \dis \frac{\varepsilon}{n + 1}$ for all $n \in \N.$ Then, it is clear that for all $n \in \N$, $A - \varepsilon^n \lel B$. Hence, for any fixed $b \in B,$ $b + \varepsilon^n \in A + C.$ Passing to the limit as $n$ goes to $+ \infty$, we get the desired conclusion. 
\end{proof}

\noindent We recall that an element $A \in \pzs$ is said to be $C$-proper if $A+C \neq Z$. All through the paper, we denote by $\ppz$ the set of such elements. We recall a fundamental result regarding $C$-proper sets.

\begin{proposition}[\cite{GLA22}]\label{propcproper}
Let $A \in \pzs$. Then, $A \in \ppz$ if and only if $A  \; \ll A.$  
\end{proposition}

\noindent Following the notation in \cite{RWE98}, we denote by $\Ns$ the set of all subsequences of $\N$ containing all $n$ beyond some given integer and by $\Nss$ the set of all subsequences of $\N.$ We now recall Painlev\'e-Kuratowski's notion of set convergence (see Khan \textit{et al.} \cite{KTZ15}). A sequence $\left\lbrace A_n \right\rbrace \subset \pzs$ is said to converge to a set $A \in \pzs$ in the sense of Painlev\'e-Kuratowski (denoted by $A_n \stackrel{PK}{\rightarrow} A$) if and only if
\begin{equation*}
{\rm Ls}(A_n)  \subset  A \subset {\rm Li}(A_n) ,
\end{equation*}
where the upper PK-limit ${\rm Ls}(A_n) $ and the lower PK-limit ${\rm Li}(A_n)$ are defined by
\begin{equation*}
{\rm Ls}(A_n)  := \limsup_{n \to \infty} A_n = \{x \in X \mid \exists N' \in \Nss,\; \exists a_n\in A_n \,(n\in N') \text{ with } a_n \stackrel{\Vert \cdot \Vert_{X}}{\longrightarrow} x\},
\end{equation*}
\begin{equation*}
\text{ and } \; {\rm Li}(A_n) := \liminf_{n \to \infty} A_n = \{x \in X \mid \exists N \in \Ns,\; \exists a_n\in A_n \,(n\in N) \text{ with } a_n \stackrel{\Vert \cdot \Vert_{X}}{\longrightarrow}  x\},
\end{equation*}
respectively. The inclusion $\dis {\rm Ls}(A_n) \subset A$ is called the upper Painlev\'e-Kuratowski convergence (denoted by $A_n \stackrel{PK}{\rightharpoonup} A$) and $A \subset {\rm Li}(A_n)$ is called the lower Painlev\'e-Kuratowski convergence (denoted by $A_n \stackrel{PK}{\rightharpoondown} A$).

\subsection{On Geoffroy spaces}

From now on, {$D$ stands for a nonempty subset of a real normed vector space $X := (X, \Vert \cdot \Vert_X)$ and $F$ represents a generic set-valued map acting from $D$ to the $C$-proper subsets of $Z$, i.e. $F: D\to \ppz.$ We point out that this last notation clearly expresses that $D \subseteq \dom F:=\{ x\in X\mid F(x)\neq\emptyset\}.$\\

\noindent Prior to undertaking any further steps, it makes sense to introduce an appropriate functional framework to conduct our work. Note that a meaningful topology on $\pzs$ when it comes to study set-valued optimization problems in the context of Kuroiwa's set approach, is the so-called \textit{set order topology}, denoted by $\topo$ (see Geoffroy and Larrouy \cite{GLA22}). As a reminder, set-valued minimization problems   are typically formulated as follows:
\begin{equation*}
{\mathcal{SOP}}(F, C, D) : \left\lbrace
\begin{array}{lc}
 {\rm Minimize~} & F(x),\\ 
 {\rm with~respect~to} & C,\\ 
{\rm and~subject~to~}  & x\in D,
\end{array}
\right.
\end{equation*}
were $F : D \to \ppz$ is a $C$-proper-valued mapping, i.e. $F(x) + C \neq Z$, for all $x \in D$. Actually, the consideration of set-valued minimization problems involving a non-$C$-proper-valued objective function is irrelevant. Since the resolution of ${\mathcal{SOP}}(F, C, D)$ rely on the comparison of the elements of the family ${\mathcal{F}_{F}}:=\{F(x)\mid\, x\in \dom F\}),$ we introduce the following conlinear topological structure for set-valued maps acting from subsets of $X$ to the subsets $\pzs$: 
\begin{equation}\label{TopoStructure}
\FF := \left(\mathcal{F}_{F}, +, \cdot, \lell, \stackrel{l}{\approx}, \ll,  \topo \right).
\end{equation}
From now on, any structure of type of \eqref{TopoStructure} will be referred to as a \textit{Geoffroy space}. In this context, the continuity of set-valued maps reads as follows.

\begin{definition}[$\topo$-continuity of set-valued maps]\label{defcont}
A set-valued map $F \in \FF$ is said to be $\topo$-continuous at $\bx\in X$ if for any $\topo$-neighborhood $\mathcal{N}_{F(\bx)}$ of $F(\bx)$ there is a neighborhood $\mathcal{U}_{\bx}$ of $\bx$ such that for all $x \in \mathcal{U}_{\bx}, \, F(x) \in \mathcal{N}_{F(\bx)}.$\\ 
One says that $F$ is $\topo$-continuous on a subset $S\subseteq X$ if it is $\topo$-continuous at any point of $S$.
\end{definition}

\begin{remark}\label{RemarkTopoCont} A set-valued map is $\topo$-continuous if and only if it is both lower and upper $\topo$-semicontinuous. These concepts are described hereinafter.
\end{remark}
%

 \begin{definition}[Lower and upper $\topo$-semicontinuity]\label{defusclsc}
 Let $F \in \FF$ be a set-valued map.\smallbreak

\noindent {\rm (a)} We say that $F$ is lower $\topo$-semicontinuous $( \topo$-lsc for short$)$ at $\bx\in X$ if for any $\varepsilon \in \inte(C)$ there is a neighborhood $\mathcal{U}_{\bx}$ of $\bx$ such that for all $x\in \mathcal{U}_{\bx}, \; F(\bx) \; \ll F(x) + \varepsilon.$\smallbreak

\noindent {\rm (b)} We say that  $F$ is upper $\topo$-semicontinuous $( \topo$-usc for short$)$ at $\bx\in X$ if 
 for any $\varepsilon \in \inte(C)$ there is a neighborhood $\mathcal{U}_{\bx}$ of $\bx$ such that for all
 $x\in \mathcal{U}_{\bx}, \; F(x) \; \ll F(\bx) + \varepsilon.$\smallbreak
\noindent We say that $F$ is $\topo$-lsc $($resp. $\topo$-usc$)$ on $S\subseteq X$ if it is $\topo$-lsc $($resp. $\topo$-usc$)$ at any point of $S$.
 \end{definition}

\noindent Here, we have gathered important properties taken from \cite[Section 4.2]{GLA22} that will be used later.

 \begin{proposition}[Characterization of lower $\topo$-semicontinuity]\label{propcharlsc}
Let $F \in \FF$ be a set-valued map and let $\bx\in X$. Then, the following assertions are equivalent:\smallbreak
  
 \noindent {\rm ~(i)} $F$ is lower $\topo$-semicontinuous at $\bx$ $;$\smallbreak

\noindent {\rm (ii)} $F$ is sequentially lower $\topo$-semicontinuous at $\bx$, i.e.,  for any sequence $\left\lbrace x_n \right\rbrace \subset X$ converging to $\bx$ and for all $\varepsilon \in \inte(C),\; F(\bx) - \varepsilon \ll F(x_n),$ for $n$ large enough. 
\end{proposition}

 \begin{proposition}[Characterization of upper $\topo$-semicontinuity]\label{propcharusc}
Let $F \in \FF$ be a set-valued map and let $\bx\in D$. Then, the following assertions are equivalent:\smallbreak
  
\noindent {\rm ~(i)} $F$ is upper $\topo$-semicontinuous at $\bx$ $;$\smallbreak

 \noindent {\rm (ii)} $F$ is sequentially upper $\topo$-semicontinuous at $\bx$, i.e., for any sequence $\left\lbrace x_n \right\rbrace \subset X$ converging to $\bx$ and for all $\varepsilon \in \inte(C), \; F(x_n) - \varepsilon \ll F(\bx),$ for $n$ large enough.
\end{proposition}

\noindent Before continuing, we introduce the concept of strong level sets which makes sense on Geoffroy spaces. 

\begin{definition}[Strong level sets]\label{defStrongLevelSet} Let $F \in \FF$ be a set-valued map and $\Omega \in \pzs.$ We call strong level set of $F$ at height $\Omega$ the set ${\rm Lev}_{\Omega}(F) := \left\lbrace x \in X \mid F(x) \; \lell \; \Omega \right\rbrace.$
\end{definition}

\begin{remark} When $F$ is $C$-closed valued Definition \ref{defStrongLevelSet} coincides with the classical concept of level set ${\rm lev}_{\Omega}(F) := \left\lbrace x \in X \mid F(x) \; \lel \; \Omega \right\rbrace.$
\end{remark}

\noindent We show that the Definition \ref{defStrongLevelSet} fits into the scalar case's classical pattern of level sets.

\begin{proposition} Let $F \in \FF$. If $F$ is lower $\topo$-semicontinous on $X$, then ${\rm Lev}_{\Omega}(F)$ is a closed subset of $X$, for all $\Omega \in \pzs.$
\end{proposition}

\begin{proof} Let $\Omega \in \pzs.$ The case ${\rm Lev}_{\Omega}(F) = \emptyset$ is trivial. Suppose that ${\rm Lev}_{\Omega}(F) \neq \emptyset.$ Let $\left\lbrace x_n \right\rbrace \subset {\rm Lev}_{\Omega}(F)$ be a sequence converging to some $\bx \in X.$ Invoking Proposition \ref{propcharlsc}, we obtain that for all $\varepsilon \in \inte(C)$, $F(\bx) - \varepsilon \; \ll F(x_n),$ eventually. Since $x_n \in {\rm Lev}_{\Omega}(F)$ for all $n \in \N$, we deduce that for all $\varepsilon \in \inte(C)$, $F(\bx) - \varepsilon \; \ll \Omega.$ Finally, applying Lemma \ref{lemmastrictStrong} $(b)$, we get the desired conclusion.
\end{proof}

\subsection{Minimality concepts}

We recall that in our context, solving ${\mathcal{SOP}}(F, C, D)$ consists in finding points $x \in D$ (called minimal solutions) such that the value of $F$ at $x$, namely $F(x)$, possesses some kind of minimality properties. Hereinafter, for any set-valued map $F \in \F$, we focus on the following concepts of minimal solutions which are of particular interest regarding real-world applications. 

\begin{definition}[Strong minimal solution] A point $\bx\in D$ is called a strong minimal solution to $\mathcal{SOP}(F, C, D)$ if for all $x \in D$, $F(\bx) \; \lel\,  F(x)$. The set of all strong minimal solutions to $\mathcal{SOP}(F, C, D)$ is defined by
\begin{equation*}
\Se(F, C,D):= \{\bx\in D \mid F(\bx) \; \lel F(x), \forall x\in D\}.
\end{equation*}
\end{definition}

\begin{definition}[Pareto minimal solution] A point $\bx\in D$ is called a Pareto minimal solution to ${\mathcal{SOP}}(F, C, D)$ if for all $x \in D$, $F(x)\lel F(\bx) \Rightarrow F(\bx)\lel F(x)$. The set of all Pareto minimal solutions to ${\mathcal{SOP}}(F, C, D)$ is defined by
\begin{center}
$\Pe(F, C, D) := \{\bx\in D \mid x \in D, F(x)\; \lel F(\bx) \Longrightarrow F(\bx)\; \lel F(x)\}.$
\end{center}
\end{definition}

\begin{definition}[Relaxed minimal solution] A point $\bx\in D$ is called a relaxed minimal solution to $\mathcal{SOP}(F, C, D)$ if there is no $x \in D$ such that $F(x)\ll F(\bx)$. The set of all relaxed minimal solutions to $\mathcal{SOP}(F, C, D)$ is the set 
\begin{equation*}
\Re(F, C,D):=\{\bx\in D \mid \nexists\, x \in D \text{ such that } F(x) \; \ll F(\bx)\}.
\end{equation*}
\end{definition}

\noindent We now introduce the concept of Geoffroy minimal solutions which is significant in Geoffroy spaces.

\begin{definition}[Geoffroy minimal solution]\label{DefGeoffroySol} A point $\bx\in D$ is called a Geoffroy minimal solution to ${\mathcal{SOP}}(F, C, D)$ if for all $x \in D$, $F(x)\lell F(\bx) \Rightarrow F(\bx)\lell F(x)$. The set of all Pareto minimal solutions to ${\mathcal{SOP}}(F, C, D)$ is defined by
\begin{center}
$\Ge(F, C, D) := \{\bx\in D \mid x \in D, F(x)\; \lell F(\bx) \Longrightarrow F(\bx)\; \lell F(x)\}.$
\end{center}
\end{definition}

\noindent We infer that:
\begin{equation}
 \Ge(F, C, D) \; \subseteq \Re(F, C, D).
\end{equation}
Indeed, assume that there is $\bx \in \Ge(F, C, D)$ such that $\bx \notin \Re(F, C, D).$ Then, there is $x_0 \in D$ such that $F(x_0) \; \ll F(\bx)$. Thus, $F(x_0) \; \lell F(\bx)$ and from Remark \ref{RemarkLinkRelations}, we get $F(x_0) \; \ll F(x_0).$ Since $F(x_0) \in \ppz$, Proposition \ref{propcproper} yields a contradiction. \\

\noindent It seems important to precise that $\Ge(F, C, D) \neq \Re(F, C, D),$ in general. 

\begin{example} Let $X = \R$, $Z = \R^2$ and $C := \R_+^2.$ We consider the map $F : \mathcal{G}_C(\R, \P_p(\R^2))$ defined by
\begin{equation*}
F(x) = \left\lbrace
\begin{array}{lcl}
\left(0, 1\right) \times \left(0, 1\right] & \text{if } & x < 0,\\ 
\left(0, 1\right) \times \left(x, 1\right] & \text{if }&  x \in \left[0, 2\right),\\ 
\left[-1,5 \right] \times \left[2, 3\right] & \text{if }&  x \geq 2.
\end{array}
\right.
\end{equation*}
Consider the unconstrained set-valued optimization problem ${\mathcal{SOP}}(F, \R_+^2, \R).$ Then, a mere computation shows that $\Ge(F, \R_+^2, \R) = \R \setminus \left(0, 2\right)$ while $\Re(F, \R_+^2, \R) = \R.$
\end{example}

\begin{remark}\label{RemInc} From Geoffroy and Larrouy \cite[\text{Proposition 3.11}]{GLA22}, we have:
\begin{equation}\label{remeffinclusion}
\text{if }  \Se(F, C,D) \neq \emptyset, \text{ then }\Se(F,C,D) = \Pe(F, C,D).
\end{equation}
In addition, when $F$ is $C$-closed-valued we have $\Pe(F,C,D) = \Ge(F,C,D).$
\end{remark}

\noindent Given Remark \ref{RemInc}, it is sufficient to restrict ourselves to the study of Geoffroy minimal solutions and relaxed minimal solutions. Other concepts being particular cases of the latter. Throughout the paper, it is assumed that $\Ge(F,C,D) \neq \emptyset$ and $\Re(F,C,D) \neq \emptyset.$

\section{Asymptotic properties of strong level sets}\label{Section 2}

Our first task is to study the upper and the lower convergence of strong level sets, in the sense of the Painlev\'e-Kuratowski convergence. Then, we will apply the obtained results to set optimization. 

\subsection{The Gamma-cone convergence}

We start by introducing a new variational convergence tailored for stability analysis in Geoffroy spaces. 

\begin{definition}[Gamma-cone convergence]\label{defgammaconv}
Consider $\bx \in X$. Let $F \in \FF$ be a set-valued map and $\left\lbrace F_n \right\rbrace \subset \pzs$ be a sequence of set-valued maps such that $F_n \in \FF,$ for all $n \in \N.$ We say that the sequence $F_n$ $\Gamma^C$-converges $($or Gamma-cone converges$)$ to $F$ at $\bx$ if and only if 
\begin{itemize}
\item[$(a)$] for all $\varepsilon \in \inte(C)$ and all sequence $\left\lbrace x_n \right\rbrace \subset X$ such that $x_n \stackrel{\Vert \cdot \Vert_{X}}{\longrightarrow} \bx$, $F(\bx)- \varepsilon \; \lell \; F_n(x_n)$ eventually.
\item[$(b)$] there is $\left\lbrace x^\star_n \right\rbrace \subset X$ with $x^\star_n \stackrel{\Vert \cdot \Vert_{X}}{\longrightarrow} \bx$ such that for all $\varepsilon \in \inte(C),\;  F_n(x^\star_n) \; \lell \; F(\bx) + \varepsilon$ eventually.
\end{itemize}
In this case, we say that $F(\bx)$ is the $\Gamma^C$-limit of $\left\lbrace F_n \right\rbrace$ at $\bx$ and we write $ \dis F_n  \underset{\bx}{\overset{\Gamma^{C}}{\longrightarrow}} F.$  When this property holds for all $\bx \in X,$ we say that $ \left\lbrace F_n \right\rbrace$ $\Gamma^C$-converges to $F$ and we write $F_n \stackrel{\Gamma^{C}}{\longrightarrow} F.$
\end{definition}

\noindent We provide an example of a Gamma-convergent sequence of set-valued maps for the sake of clarity.

\begin{example} Let $X = \R$, $Z = \R$ and $C := \R_+.$ We consider the map $F \in \mathcal{G}_C(\R, \P_p(\R))$ defined by
\begin{equation*}
F(x) = \left(\cos(x), 1 + \cos(3x) + \cos(5x) \right). 
\end{equation*}
We infer that the sequence $\left\lbrace F_n \right\rbrace$ of set-valued maps given by 
\begin{equation*}
F_n(x) = \left\lbrace \cos \left(x + \frac{1}{n+1}\right) \right\rbrace, \; \text{ for all } n \in \N,  
\end{equation*}
is $\Gamma^{\R_{+}}$-convergent to $F$ at $0.$ Indeed, let $\varepsilon \in \inte(\R_{+}).$ Assume that $\left\lbrace x_n \right\rbrace \subset \R$ is such that $x_n \stackrel{\vert \cdot \vert}{\longrightarrow} 0$. Then, from the continuity of $\cos(\cdot)$ on $\R$ we have $\dis \cos \left(x_n + \frac{1}{n+1}\right) \in \left(\cos(0) - \varepsilon, \cos(0) + \varepsilon \right),$ eventually. Therefore, it follows that 
\begin{equation*}
F_n(x_n) = \left\lbrace \cos \left(x_n + \frac{1}{n+1}\right) \right\rbrace \; \subset \left[ \cos(0) - \varepsilon, + \infty \right) =: {\rm cl}(F(0) - \varepsilon + \R_{+}), \text{ eventually.}
\end{equation*}
In other words, $F(0) - \varepsilon \; \lell F_n(x_n),$ eventually. Thus, property $(a)$ in previous Definition \ref{defgammaconv} is satisfied. Now, consider the sequence $\left\lbrace x_n^\star \right\rbrace \subset \R$ given by $x_n^\star = \frac{1}{n+ 1},$ for all $n \in \N.$ Then, $x_n^\star \stackrel{\vert \cdot \vert}{\longrightarrow} 0$ and since $\cos(x_n^\star + \frac{1}{n+ 1}) \; \leq \cos(0),$ for all $n \in \N,$ we get 
\begin{equation*}
F(0) + \varepsilon = \left(1 + \varepsilon, 3 + \varepsilon \right) \; \subset \left[ \cos(x_n^\star + \frac{1}{n+ 1}), + \infty \right) =: {\rm cl}(F_n(x_n^\star)+ \R_{+}), \text{ eventually.}
\end{equation*}
In other words, $F_n(x^\star_n) \; \lell \; F(\bx) + \varepsilon$ eventually. Hence, property $(b)$ in Definition \ref{defgammaconv} is satisfied. Finally, we have $ \dis F_n  \underset{0}{\overset{\Gamma^{\R_{+}}}{\longrightarrow}} F.$
\end{example}

\begin{remark} The properties $(a)$ and $(b)$ in Definition \ref{defgammaconv} will be referred to as the lower and the upper Gamma-convergence of $F_n$ to $F$ at $\bx$, respectively.
\end{remark}

\noindent As we now show, the lower $\Gamma^C$-convergence can be characterized in terms of neighborhoods. 

\begin{lemma}[Characterization of the lower Gamma-cone convergence]\label{lemmaTopoGammaConv} Let $F \in \FF$ and $\left\lbrace F_n \right\rbrace \subset \pzs$ be such that $F_n \in \FF,$ for all $n \in \N.$ Assume that $F_n  \underset{\bx}{\overset{\Gamma^{C}}{\longrightarrow}} F$ for some $\bx \in X.$ Then, the following assertions are equivalent.
\begin{itemize}
\item[$(i)$] Condition $(a)$ in Definition \ref{defgammaconv} holds for $\bx.$ 
\item[$(ii)$] For all $\varepsilon \in \inte(C),$ there exists a neighborhood $\mathcal{N}_{\varepsilon, \bx}$ of $\bx$ $($depending on $\varepsilon )$ such that for all $x_0 \in \mathcal{N}_{\varepsilon, \bx},$ $F(\bx) - \varepsilon \; \ll F_n(x_0)$ eventually.
\end{itemize}
\end{lemma}

\begin{proof} Let us prove that $(i)$ implies $(ii)$. Assume to the contrary that $(ii)$ is false. Hence, we can find $\varepsilon_0 \in \inte(C)$ such that we can construct a sequence $\left\lbrace x_n \right\rbrace \subset X$ such that $x_n \in \ball_{\frac{1}{n}}(\bx)$ for all $n \in \N$ and $F(\bx) - \varepsilon_0 \; \nlell F_n(x_n),$ i.e. it is false that $F(\bx) - \varepsilon_0 \; \lell F_n(x_n)$. Since $x_n \stackrel{\Vert \cdot \Vert_{X}}{\longrightarrow} \bx$, the proof is complete. Now, suppose that $(ii)$ holds and take $\varepsilon \in \inte(C)$. Then, we can find a neighborhood $\mathcal{N}_{\varepsilon, \bx}$ of $\bx$ such that  for all $x_0 \in \mathcal{N}_{\varepsilon, \bx},$ $F(\bx) - \varepsilon \; \ll F_n(x_0),$ eventually. If $\left\lbrace x_n \right\rbrace \subset X$ is a sequence converging to $\bx$ then $x_n \in \mathcal{N}_{\varepsilon, \bx}$ for $n$ large enough. Therefore, we have $F(\bx) - \varepsilon \; \lell F_n(x_n),$  eventually.
\end{proof}

\noindent It is worth mentioning that the $\Gamma^C$-limit is not unique. More precisely, we have the following result.

\begin{proposition} Let $F, G \in \FF$ and let $\left\lbrace F_n \right\rbrace \subset \pzs$ be such that $F_n \in \FF,$ for all $n \in \N.$ Assume that $F_n  \underset{\bx}{\overset{\Gamma^{C}}{\longrightarrow}} F$ for some $\bx \in X.$ Then, $F_n  \underset{\bx}{\overset{\Gamma^{C}}{\longrightarrow}} G$ if and only if $F(\bx) \stackrel{l}{\approx} G(\bx).$ 
\end{proposition}

\begin{proof} If $F(\bx) \stackrel{l}{\approx} G(\bx),$ it is clear that $F_n  \underset{\bx}{\overset{\Gamma^{C}}{\longrightarrow}} G$. Suppose now that $F_n  \underset{\bx}{\overset{\Gamma^{C}}{\longrightarrow}} G$ and take $\varepsilon \in \inte(C).$ Then, there exists a sequence $\left\lbrace x^\star_n \right\rbrace \subset X$ converging to $\bx$ and such that $F_n(x_n^\star) \; \lell G(\bx) + \dis \frac{\varepsilon}{4},$ eventually. From Lemma \ref{lemmastrictStrong} $(a)$, we have
\begin{equation}\label{EqPropNonUniqueness1}
F_n(x_n^\star) \; \ll G(\bx) + \frac{\varepsilon}{2}, \; \text{ eventually}.
\end{equation}
Now, since $F_n  \underset{\bx}{\overset{\Gamma^{C}}{\longrightarrow}} F$, we deduce by Lemma \ref{lemmastrictStrong} $(a)$ again that $\dis F(\bx) - \frac{\varepsilon}{2} \; \ll F_n(x_n^\star),$ eventually. This, together with \eqref{EqPropNonUniqueness1} yields that $F(\bx) - \varepsilon \, \ll G(\bx),$ for all $\varepsilon \in \inte(C).$  Thus, from Lemma \ref{lemmastrictStrong} $(b)$, we obtain that $F(\bx) \; \lell G(\bx).$ By permuting $F$ and $G$ throughout the proof, we prove that we also have $G(\bx) \; \lell F(\bx).$.
\end{proof}

\noindent As in the scalar case, $\Gamma^C$-limits are always lower $\topo$-semicontinuous. 

\begin{theorem}[Lower $\topo$-semicontinuity of the $\Gamma^C$-limits] Let $F \in \FF$ be a set-valued map and $\left\lbrace F_n \right\rbrace \subset \pzs$ be a sequence of set-valued maps with $F_n \in \FF$ for all $n \in \N$ and such that $F_n  \overset{\Gamma^{C}}{\longrightarrow} F.$ Then, $F$ is lower $\topo$-semicontinuous.
\end{theorem}

\begin{proof} Take $\varepsilon \in \inte(C)$ and $\bx \in X$. From Lemma \ref{lemmaTopoGammaConv}, there exists $\delta_\varepsilon > 0$ such that for all $x_0 \in \ball_{\delta_{\varepsilon}}(\bx)$:
\begin{equation}\label{IneqGammaSemicontinuous1}
F(\bx) - \frac{\varepsilon}{2} \; \ll F_n(x_0), \; \text{eventually}.
\end{equation}
Let us fix an element $\hx \in \ball_{\delta_{\varepsilon}}(\bx).$ Then, there exists a sequence $\left\lbrace z_n \right\rbrace \subseteq X$ converging to $\hx$ such that 
\begin{equation}\label{IneqGammaSemicontinuous2}
F_n(z_n) \; \lell F(\hx) + \frac{\varepsilon}{2}, \; \text{eventually}.
\end{equation}
Furthermore, we can find a pair $(\beta, n_0) \in \R_+^\star \times \N$ such that $z_{n_{0}} \in \ball_{\beta}(\hx)$ and $\ball_{\beta}(\hx) \subset \ball_{\delta_{\varepsilon}}(\bx).$ Hence, combining \eqref{IneqGammaSemicontinuous1} and \eqref{IneqGammaSemicontinuous2}, we deduce that for all $x \in \ball_{\delta_{\varepsilon}}(\bx)$, $F(\bx) - \varepsilon \; \ll F(x).$ 
\end{proof}

\noindent In addition, we prove that the $C$-lower stationary sequence defined by $F_n \stackrel{l}{\approx} F,$ $\Gamma^C$-converges to $F$ when $F$ is $\topo$-continuous.

\begin{theorem} Let $F \in \FF$ be a set-valued map and $\left\lbrace F_n \right\rbrace \subset \pzs$ be a sequence of set-valued maps with $F_n \stackrel{l}{\approx} F,$ for all $n \in \N.$ If $F$ is $\topo$-continuous at some $\bx \in X$, then $\dis F_n  \underset{\bx}{\overset{\Gamma^{C}}{\longrightarrow}} F.$
\end{theorem}

\begin{proof} Assume that $F$ is $\topo$-continuous at $\bx \in X.$  According to Remark \ref{RemarkTopoCont}, the set-valued map $F$ is both lower $\topo$-semicontinuous and upper $\topo$-semicontinuous at $\bx$. Since $F(\bz) \; \lell F_n(\bz)$ for all $\bz \in X,$ we deduce from Proposition \ref{propcharlsc}, that for all $\varepsilon \in \inte(C)$ and all sequence $\left\lbrace x_n \right\rbrace \subset X$ such that $x_n \stackrel{\Vert \cdot \Vert_{X}}{\longrightarrow} \bx$, $F(\bx)- \varepsilon \lell F_n(x_n)$ eventually. On the other hand, we can consider a sequence $\left\lbrace z_n \right\rbrace \subset X$ such that $z_n \stackrel{\Vert \cdot \Vert_{X}}{\longrightarrow} \bx.$ Because $F_n(\by) \; \lell F(\by)$ for all $\by \in X,$ we deduce from Proposition \ref{propcharusc} that for all $\varepsilon \in \inte(C),$ $F_n(z_n) \; \lell \; F(\bx) + \varepsilon$ eventually. The proof is complete.
\end{proof}

\subsection{Upper and lower convergence of strong level sets}

\noindent Thanks to our new variational convergence for sequence of set-valued maps, we are so far able to prove the main results of this section. First, we study the upper convergence of strong level sets in Geoffroy spaces.

\begin{theorem}[Upper convergence of strong level sets]\label{thstablevUpper}
Let $F \in \FF$ be a set-valued map and $\Omega\in\pzs.$ Consider sequences $\left\lbrace F_n \right\rbrace, \left\lbrace \Omega_n \right\rbrace \subset \pzs$ such that $F_n \in \FF$ and $\Omega_n \in \pzs$ for all $n \in \N.$ Assume that
\begin{itemize}
\item[$(a)$] $F_n \stackrel{\Gamma^{C}}{\longrightarrow} F$ $;$
\item[$(b)$] for all $\varepsilon \in \inte(C),$ there is $n_k \in \Nss$ such that $\Omega_{n_{k}} - \varepsilon \; \lell \Omega$ eventually.
\end{itemize}
Then, $\dis {\rm Lev}_{\Omega_{n}}(F_{n}) \stackrel{PK}{\rightharpoonup} \dis {\rm Lev}_{\Omega}(F).$
\end{theorem}

\begin{proof} Let $\bx \in {\rm Ls}\left({\rm Lev}_{\Omega_{n}}(F_{n})\right)$. Then, we can find a sequence $n_k \in \Nss$ such that $x_{n_{k}} \stackrel{\Vert \cdot \Vert_{X}}{\longrightarrow} \bx$ with $\xnk\in {\rm Lev}_{\Omega_{n_k}}(F_{n_k}),$ for all $k \in \N$. Thus, it follows that
\begin{equation}\label{eqsuplevT1}
F_{n_k}(\xnk) \; \lell \Omega_{n_k}, \text{ for all } k\in\N.
\end{equation}   
Let $\varepsilon \in\inte(C)$. Since $\left\lbrace F_n \right\rbrace$ $\Gamma^C$-converges to $F$ at $\bx$, we deduce that 
\begin{equation}\label{eqsuplevT2}
F(\bx) - \dis \frac{\varepsilon}{4} \; \lell F_{n_k}(\xnk), \; \text{  eventually.}
\end{equation}
Combining $(b)$ with \eqref{eqsuplevT1}, we know that there is $n_{k_{l}} \in \Nss$ such that $F_{n_{k_{l}}}(x_{n_{k_{l}}}) - \dis \frac{\varepsilon}{4} \; \lell \Omega,$ eventually. This, together with \eqref{eqsuplevT2} gives that $F(\bx) - \dis \frac{\varepsilon}{2} \; \lell \Omega.$ Invoking Lemma \ref{lemmastrictStrong} $(a)$ we get $F(\bx) - \varepsilon \; \ll \Omega.$ As $\varepsilon$ has been chosen arbitrarily, we infer from the Lemma \ref{lemmastrictStrong} $(b)$ that $F(\bx) \; \lell \Omega,$ i.e. $x\in {\rm Lev}_\Omega(F).$
\end{proof}

\noindent One can see that the concept of strong level sets is closely connected to that of Geoffroy minimal solution. Indeed, we have the following result.

\begin{proposition}\label{PropGeof} Let $\bx \in \Ge(F,C,D).$ Then, ${\rm Lev}_{F(\bx)}(F) \subseteq \Ge(F,C,D).$
\end{proposition}

\begin{proof} Let $\bx \in \Ge(F,C,D).$ Assume that $\bz \in D$ is such that $F(\bz) \; \lell F(x_0)$ with $x_0 \in {\rm Lev}_{F(\bx)}(F).$ Then, $F(\bz) \; \lell F(\bx).$ Since $\bx$ is a Geoffroy minimal solution, we have $F(\bx) \; \lell F(\bz).$ Therefore, the very nature of $x_0$ yields that $x_0 \in \Ge(F,C,D).$
\end{proof}

\noindent As a result, we obtained the following stability result for the set of Geoffroy minimal solutions.

\begin{theorem}[External stability of Geoffroy minimal solutions set]\label{ExternStabGeof} Consider a sequence $\left\lbrace F_n \right\rbrace \subset \pzs$ such that $F_n \in \F,$ for all $n \in \N.$ Assume that:
\begin{itemize}
\item[$(a)$] $F_n \stackrel{\Gamma^{C}}{\longrightarrow} F$ $;$
\item[$(b)$] $F$ is upper $\topo$-semicontinuous at some $\bx \in \Ge(F,C,D).$
\end{itemize}
Then, $\dis {\rm Lev}_{F(x_{n})}(F_{n}) \stackrel{PK}{\rightharpoonup} \Ge(F,C,D)$ for all sequence $\left\lbrace x_n \right\rbrace \subset X$ converging to $x_0 \in {\rm Lev}_{F(\bx)}(F).$
\end{theorem}

\begin{proof} Let $x_0 \in {\rm Lev}_{F(\bx)}(F)$ with $\bx \in \Ge(F,C,D).$ Since $F$ is upper $\topo$-semicontinuous at some $\bx \in \Ge(F,C,D)$, it is sequentially upper $\topo$-semicontinuous at $x_0$ (see Proposition \ref{propcharusc}). Therefore, if $\left\lbrace x_n \right\rbrace \subset X$ is such that $x_n \stackrel{\Vert \cdot \Vert_{X}}{\longrightarrow} x_0,$ it follows that for all $\varepsilon \in \inte(C),$ $F(x_n) - \varepsilon \lell F(x_0),$ eventually. Thus, invoking Theorem \ref{thstablevUpper}, we deduce that $\dis {\rm Lev}_{F(x_{n})}(F_{n})  \stackrel{PK}{\rightharpoonup} {\rm Lev}_{F(x_{0})}(F)$.  Hence, it follows from the identity \eqref{Rep1} that $\dis {\rm Lev}_{F(x_{n})}(F_{n}) \stackrel{PK}{\rightharpoonup} \Ge(F,C,D).$ Applying Proposition \ref{PropGeof}, the proof is complete.
\end{proof}

\noindent Now, we focus on the lower convergence of strong level sets.

\begin{theorem}[Lower convergence of strong level sets]\label{thstablevLower}
Let $F \in \FF$ be a set-valued map and $\Omega\in\pzs.$ Consider sequences $\left\lbrace F_n \right\rbrace, \left\lbrace \Omega_n \right\rbrace \subset \pzs$ such that $F_n \in \FF$ and $\Omega_n \in \pzs$ for all $n \in \N.$ Assume that
\begin{itemize}
\item[$(a)$] $F_n \stackrel{\Gamma^{C}}{\longrightarrow} F$ $;$
\item[$(b)$] $\Omega \; \ll \Omega_{n},$ eventually.
\end{itemize}
Then, $\dis {\rm Lev}_{\Omega_{n}}(F_{n}) \stackrel{PK}{\rightharpoondown} {\rm Lev}_{\Omega}(F).$
\end{theorem}
\begin{proof} First, let us observe that from $(b)$, there is $\varepsilon_0 \in \inte(C)$ such that 
\begin{equation}\label{LowerStabEq1}
\Omega + \varepsilon_0 \; \lell \Omega_{n}, \text{ eventually.}
\end{equation}
Let $\bx \in {\rm Lev}_{\Omega}(F)$. According to $(a)$, there is a sequence $\left\lbrace x_n \right\rbrace \subset X$ converging to $\bx$ such that $F_n(x_n) \; \lell \; F(\bx) + \varepsilon_0,$ eventually. Because $\bx \in {\rm Lev}_{\Omega}(F)$, it follows that $\dis F_n(x_n) \; \lell \; \Omega + \varepsilon_0,$ eventually. Combining the latter inequality with \eqref{LowerStabEq1}, we deduce that for $n$ large enough, $F_n(x_n) \; \lell \; \Omega_n.$ 
\end{proof}

\noindent To study prove the internal stability of the set of Geoffroy minimal solutions, we introduce the following notion. In what follows, we write $\Te(F,C,D)$ to indicate the set of minimal solutions to the problem ${\mathcal{SOP}}(F, C, D)$ of type $\mathcal{T}$ where $\mathcal{T} \in \left\lbrace \mathcal{G}, \mathcal{P}, \mathcal{S}, \mathcal{R} \right\rbrace.$

\begin{definition} Let $\mathcal{T} \in \left\lbrace \mathcal{G}, \mathcal{P}, \mathcal{S}, \mathcal{R} \right\rbrace.$ We say that $\Te(F,C,D)$ admits a countable number of representant if and only if there exists $x^1, x^2, \ldots, x^n \in D$ $(n \in \N)$ such that 
\begin{equation}\label{Rep1}
\Te(F,C,D) \; = \; {\rm Lev}_{F(x^{1})}(F) \; \cup {\rm Lev}_{F(x^{2})}(F) \; \cup \; \ldots \; \cup {\rm Lev}_{F(x^{n})}(F).
\end{equation} 
with 
\begin{equation}\label{Rep2}
{\rm Lev}_{F(x^{i})}(F) \; \cap \; {\rm Lev}_{F(x^{j})}(F) = \emptyset, \; \text{ for all } i \neq j.
\end{equation} 
When \eqref{Rep1} and \eqref{Rep2} hold, the vector $\left(x^1, x^2, \ldots, x^n \right) \in D^n$ is called a $\T$-representant of $\Te(F,C,D).$ 
\end{definition} 

\noindent Therefore, we were able to prove the next theorem. 

\begin{theorem}[Internal stability of Geoffroy minimal solutions set]\label{InternStabGeof} Assume that $\Ge(F,C,D)$ has a countable number of representant and let $\left(x^1, x^2, \ldots, x^n \right) \in D^n$ be one of its $\mathcal{G}$-representant. Consider a sequence $\left\lbrace F_n \right\rbrace \subset \pzs$ such that $F_n \in \F,$ for all $n \in \N.$ Assume that
\begin{itemize}
\item[$(a)$] $F_n \stackrel{\Gamma^{C}}{\longrightarrow} F$ $;$
\item[$(b)$] there exists a sequence $\left\lbrace x_n \right\rbrace \subset X$ such that $F(x^i) \; \ll F(x_n)$ eventually, for all $i \in \left\lbrace 1, 2, \ldots, n \right\rbrace.$
\end{itemize}
Then, $\dis {\rm Lev}_{F(x_{n})}(F_{n}) \stackrel{PK}{\rightharpoondown} \Ge(F,C,D).$
\end{theorem}

\begin{proof} From Theorem \ref{thstablevLower}, $\dis {\rm Lev}_{F(x_{n})}(F_{n})  \stackrel{PK}{\rightharpoondown} {\rm Lev}_{F(x^i)}(F)$ for each component $x^i$ of the $\mathcal{G}$-representant of $\Ge(F,C,D).$ Hence, it follows from the identity \eqref{Rep1} that $\dis {\rm Lev}_{F(x_{n})}(F_{n}) \stackrel{PK}{\rightharpoondown} \Ge(F,C,D).$
\end{proof}

\begin{remark} Obviously, from Remark \ref{RemInc}, Theorems \ref{ExternStabGeof} and \ref{InternStabGeof} are valid for 
$\Se(F,C,D)$ and $\Pe(F,C,D)$ when $F$ is $C$-closed valued. 
\end{remark}

\section{Stability results for set-valued minimization problems}\label{Section 3}

Let us now consider the sequence $\left\lbrace {\mathcal{SOP}}(F_n, C, D_n) \right\rbrace$ of set-valued optimization problems given by
\begin{equation*}
{\mathcal{SOP}}(F_n, C, D_n) : \left\lbrace
\begin{array}{lc}
 {\rm Minimize~} & F_n(x),\\ 
 {\rm with~respect~to} & C,\\ 
{\rm and~subject~to~}  & x\in D_n,
\end{array}
\right.
\end{equation*}
where $F_n \in \DNFF,$ for all $n \in \N.$ In the remainder of the paper, our objective is to define suitable approximations to the minimal solutions to ${\mathcal{SOP}}(F, C, D)$ by considering a sequence of approximate set-valued optimization problems of the form ${\mathcal{SOP}}(F_n, C, D_n)$. Indeed, we are interested in stability analysis for set-valued optimization problems in the sense of Guti\'errez \textit{et al.} \cite{GMM16}.  More precisely, we study the external and internal stability of minimal solutions to ${\mathcal{SOP}}(F, C, D).$ External stability refers to the fact that the limit of a convergent sequence of minimal solutions to ${\mathcal{SOP}}(F_n, C, D_n)$ is a solution to ${\mathcal{SOP}}(F, C, D)$, while internal stability refers to the fact that a given solution to ${\mathcal{SOP}}(F, C, D)$ can be expressed as a limit of solutions to ${\mathcal{SOP}}(F_n, C, D_n).$ It is worth mentioning that in our work we consider the perturbation both on the objective function and on the admissible domain. \\

To carry out our work, we need to extend the Definition \ref{defgammaconv} to a notion of Gamma-cone convergence for sequences of set-valued maps whose domain depends on the index of each element. For this purpose, we shall first introduce the concept of Kuratowski pair. 

\begin{definition}[Kuratowski pair] Let $\left\lbrace D_n \right\rbrace, D \subseteq X$. We say that the pair $\left( \left\lbrace D_n \right\rbrace, D \right)$ is a Kuratowski pair if and only if for any sequence $\left\lbrace x_n \right\rbrace \subset X$ with $x_n \in D_n$ for all $n \in \N,$ it is possible to find a pair $(n_k, \bx) \in \Nss \times D$ such that $x_{n_{k}} \stackrel{\Vert \cdot \Vert_{X}}{\longrightarrow} \bx.$
\end{definition}

\noindent We provide a useful criterion to Kuratowski pair.

\begin{proposition}\label{criterionKuratowskiPair} Let $\left\lbrace D_n \right\rbrace, D \subseteq X$. Assume that $D_n \stackrel{PK}{\rightharpoonup} D$. If there exists $(n_0, S) \in \N \times \pxs$ such that $D_n  \subseteq S$ with ${\rm cl}(S)$ compact for all $n \geq n_0$, then $\left( \left\lbrace D_n \right\rbrace, D \right)$ is a Kuratowski pair.
\end{proposition}

\begin{proof} Let $\left\lbrace x_n \right\rbrace \subset X$ with $x_n \in D_n,$ for all $n \in \N.$ By assumption, we have $\left\lbrace x_n \right\rbrace \subset {\rm cl}(S)$, for $n$ large enough. Moreover, since ${\rm cl}(S)$ is compact, there is a pair $(n_k, \bx) \in \Nss \times {\rm cl}(S)$ such that $x_{n_{k}} \stackrel{\Vert \cdot \Vert_{X}}{\longrightarrow} \bx$ and then $\bx \in \dis {\rm Ls}(D_n).$ Invoking the fact that $D_n \stackrel{PK}{\rightharpoonup} D$, we conclude that $\bx \in D.$
\end{proof}

\noindent Now, we are able to introduce a more general Gamma-cone convergence in Geoffroy spaces. 

\begin{definition}[Sequential Gamma-cone convergence]\label{defgammaconv}
Assume that $\left\lbrace D_n \right\rbrace, D \subseteq X$ and $\bx \in D$.\\ 
Let $F \in \FF$ and $\left\lbrace F_n \right\rbrace \subset \ppz$ be a sequence of set-valued maps such that $ F_n  \in \FF$ for all $n \in \N$. We say that $\left\lbrace F_n \right\rbrace$ $\Gamma^C_{\text{seq}}$-converges $($or sequentially Gamma-cone converges$)$ to $F$ at $\bx$ if and only if 
\begin{itemize}
\item[$(a)$] $\left( \left\lbrace D_n \right\rbrace, D \right)$ is a Kuratowski pair $;$
\item[$(b)$] for all $\varepsilon \in \inte(C)$ and all $\left\lbrace x_n \right\rbrace \subset X$ with $x_n \in D_n$ for all $n \in \N$ satisfying $x_n \stackrel{\Vert \cdot \Vert_{X}}{\longrightarrow} \bx$, we have:
\begin{equation*}
F(\bx)- \varepsilon \; \lell \; F_n(x_n), \text{ eventually } ;
\end{equation*}
\item[$(c)$] there exists $\left\lbrace x_n^\star \right\rbrace \subset X$ with $x_n^\star \in D_n$ for all $n \in \N$ converging to $\bx$ such that for all $\varepsilon \in \inte(C): $ 
\begin{equation*}
F_n(x^\star_n) \; \lell \; F(\bx) + \varepsilon,\;  \text{ eventually.}
\end{equation*}
\end{itemize}
In this case, we write $ \dis F_n  \underset{\bx}{\overset{\Gamma^{C}_{\text{seq}}}{\longrightarrow}} F.$ When this property holds for all $\bx \in D,$ we write $F_n \underset{D_n, D}{\overset{\Gamma^{C}_{\text{seq}}}{\longrightarrow}} F.$
\end{definition}

\noindent We start studying the external stability of minimal solutions to ${\mathcal{SOP}}(F, C, D)$. As we now show, the sequential Gamma-cone convergence is enough to establish the external stability of relaxed minimal solution to ${\mathcal{SOP}}(F, C, D).$

\begin{theorem}[External stability of relaxed minimal solutions]\label{UpperStabilityRelaxed} Consider the set-valued optimization problem ${\mathcal{SOP}}(F, C, D)$ and let $\left\lbrace {\mathcal{SOP}}(F_n, C, D_n) \right\rbrace$ be a sequence of approximated set-valued optimization problems. If $F_n \underset{D_n, D}{\overset{\Gamma^{C}_{\text{seq}}}{\longrightarrow}} F$, then for all sequence $\dis \left\lbrace x_n \right\rbrace \subset X$ such that $x_n \in \Re(F_n, C, D_n)$ for all $n \in \N$, there exists $(n_k, \bx) \in \Nss \times \Re(F, C, D)$ such that $x_{n_{k}} \stackrel{\Vert \cdot \Vert_{X}}{\longrightarrow} \bx.$
\end{theorem}

\begin{proof} Let $\dis \left\lbrace x_n \right\rbrace \subset X$ be a sequence such that $x_n \in   \Re(F_n, C, D_n)$, for all $n \in \N$. Then, $x_n \in D_n$ for all $n\in \N$. Since $F_n \underset{D_n, D}{\overset{\Gamma^{C}_{\text{seq}}}{\longrightarrow}} F$, $\left( \left\lbrace D_n \right\rbrace, D \right)$ is a Kuratowski pair. Thus, there exists $\left( n_k; \bx \right) \in \Nss \times D$ such that  $x_{n_{k}} \stackrel{\Vert \cdot \Vert_{X}}{\longrightarrow} \bx.$ Let us show that $\bx \in \Re(F, C, D).$ Assume to the contrary that there is $x_0 \in D$ such that $F(x_0) \; \ll F(\bx)$. Then, we can find $\varepsilon_0 \in \inte(C)$ such that 
\begin{equation}\label{EqRmin1}
F(x_0) + \varepsilon_0 \; \lel F(\bx).
\end{equation}
From now on, we take $\varepsilon, \varepsilon^\prime \in \inte(C)$ satisfying $\varepsilon_0 = \varepsilon + \varepsilon^\prime.$ Because $x_{n_{k}} \stackrel{\Vert \cdot \Vert_{X}}{\longrightarrow} \bx$ and $F_n \underset{D_n, D}{\overset{\Gamma^{C}_{\text{seq}}}{\longrightarrow}} F,$ we get
\begin{equation}\label{EqRmin2}
F(\bx) - \varepsilon \; \ll F_{n_{k}}(x_{n_{k}}), \; \text{ eventually.}
\end{equation}
Now, using again the $\Gamma^C_{\text{seq}}$-convergence of $\left\lbrace F_n \right\rbrace$ to $F$, we know that there exists a sequence $\left\lbrace \varphi_n \right\rbrace \subseteq X$ with $\varphi_n \in D_n$ for all $n \in \N$ converging to $x_0$ such that  
\begin{equation}\label{EqRmin3}
F_{n_{k}}(\varphi_{n_{k}}) \; \ll F(x_0) + \varepsilon^\prime, \; \text{ eventually.}
\end{equation}
Therefore, combining \eqref{EqRmin1}, \eqref{EqRmin2} and \eqref{EqRmin3} with the choosen decomposition for $\varepsilon_0$, we deduce that
\begin{align*}
F_{n_{k}}(\varphi_{n_{k}}) \; \ll F(x_0) + \varepsilon^\prime \lel \; F(\bx) - \varepsilon \; \ll  F_{n_{k}}(x_{n_{k}}), \; \text{ eventually.}
\end{align*}
Unfortunately, this contradicts the very nature of the sequence $\left\lbrace x_n \right\rbrace$.
\end{proof}

\noindent Now, we turn our attention to the case of Geoffroy minimal solutions. To this end, we introduce a new tool, which we refer to as the sequential lower converse property.  

\begin{definition}[Sequential lower converse property] Assume that $\left\lbrace D_n \right\rbrace, D \subseteq X$. Let $F \in \FF$ and let $\left\lbrace F_n \right\rbrace \subset \pzs$ be a sequence of set-valued maps such that $F_n  \in \FF$ for all $n \in \N$. We say that the pair $\left( (\left\lbrace F_n \right\rbrace, \left\lbrace D_n \right\rbrace), (F, D) \right)$ has the sequential lower converse property if and only if for any sequences $\left\lbrace x_n \right\rbrace, \left\lbrace \varphi_n \right\rbrace \subseteq X$ such that $x_n, \varphi_n \in D_n$ for all $n \in \N,$ having $x_{n} \stackrel{\Vert \cdot \Vert_{X}}{\longrightarrow} \bx$ and $\varphi_{n} \stackrel{\Vert \cdot \Vert_{X}}{\longrightarrow} x_0$ for some $\bx, x_0 \in D$ with $F(\bx) \; \lell F(x_0)$ implies that $F_n(x_n) \; \lell F_n(\varphi_n), \text{ eventually}.$
\end{definition}

\noindent Using this notion, we obtain the following stability result for Geoffroy minimal solutions.

\begin{theorem}[External stability of Geoffroy minimal solutions]\label{TheoUppGeofStab} Consider the set-valued optimization problem ${\mathcal{SOP}}(F, C, D)$ and let $\left\lbrace {\mathcal{SOP}}(F_n, C, D_n) \right\rbrace$ be a sequence of approximated set-valued optimization problems. Assume that
\begin{enumerate}
\item[$(a)$] $F_n \underset{D_n, D}{\overset{\Gamma^{C}_{\text{seq}}}{\longrightarrow}} F$ $;$
\item[$(b)$] the pair $\left( (\left\lbrace F_n \right\rbrace, \left\lbrace D_n \right\rbrace), (F, D) \right)$ has the sequential lower converse property.
\end{enumerate}
Then, for all sequence $\dis \left\lbrace x_n \right\rbrace \subset X$ satisfying $x_n \in   \Ge(F_n, C, D_n)$ for all $n \in \N$, there exists a pair $(n_k, \bx) \in \Nss \times \Ge(F, C, D)$ such that $x_{n_{k}} \stackrel{\Vert \cdot \Vert_{X}}{\longrightarrow} \bx.$
\end{theorem}

\begin{proof} Let $\dis \left\lbrace x_n \right\rbrace \subset X$ be a sequence satisfying $x_n \in \Ge(F_n, C, D_n)$, for all $n \in \N$. Proceeding as in the proof of Theorem \ref{UpperStabilityRelaxed}, we prove that it is possible to find a pair $(n_k, \bx) \in \Nss \times D$ such that $x_{n_{k}} \stackrel{\Vert \cdot \Vert_{X}}{\longrightarrow} \bx.$ Let us show that $\bx \in \Ge(F, C, D).$ Let $\bz \in D$ be such that $F(\bz) \; \lell F(\bx).$ We fix $\varepsilon \in \inte(C)$. Since $F_n \underset{D_n, D}{\overset{\Gamma^{C}_{\text{seq}}}{\longrightarrow}} F$, we can find a sequence $\left\lbrace z_n \right\rbrace \subseteq X$ with $z_n \in D_n$ for all $n \in \N$ converging to $\bz$ and such that 
\begin{equation}\label{EqParetoUpperStab2}
F_n(z_n) \; \lell F(\bz) + \frac{\varepsilon}{4}, \; \text{ eventually.}
\end{equation}
From assumption $(b)$, we deduce that $F_n(z_n) \; \lell F_n(x_n),$ eventually. Because $x_n \in \Ge(F_n, C, D_n)$, for all $n \in \N,$ we get $F_n(x_n) \; \lell F_n(z_n),$ eventually. Now, since $x_{n_{k}} \stackrel{\Vert \cdot \Vert_{X}}{\longrightarrow} \bx,$ it follows from $(a)$ that 
\begin{equation}\label{EqParetoUpperStab3}
F(\bx) - \frac{\varepsilon}{4}\; \lell F_{n_{k}}(x_{n_{k}}), \; \text{ eventually.}
\end{equation}
Therefore, combining \eqref{EqParetoUpperStab2} and \eqref{EqParetoUpperStab3}, we obtain $F(\bx) - \varepsilon \; \ll F(\bz).$ Because $\varepsilon$ has been chosen arbitrarily, we conclude from the Lemma \ref{lemmastrictStrong} $(b)$ that $F(\bx) \; \lell F(\bz).$ In other words, $\bx \in \Ge(F,C,D).$
\end{proof}

In what follows, we focus on the internal stability of minimal solutions to ${\mathcal{SOP}}(F, C, D)$. Before all, we formulate the following hypothesis:

\begin{assumption}{H}\label{Assumption A}
Let $\bx \in D$. If $\Te(F,C,D) \neq \emptyset$, then there exists $x^\star \in \lev_{F(\bx)}(F)$ such that $x^\star \in \Te(F,C,D).$ 
\end{assumption}

\noindent Now, we show that any Geoffroy minimal solution to ${\mathcal{SOP}}(F, C, D)$ can be approximated by a sequence of Geoffroy minimal solution to the approximated problems ${\mathcal{SOP}}(F_n, C, D_n).$ To establish our result, we combine the sequential $\Gamma^C$-convergence with the aforementioned domination property.

\begin{theorem}[Internal stability of Geoffroy minimal solutions]\label{TheoLowGeofStab} Consider the set-valued optimization problem ${\mathcal{SOP}}(F, C, D)$ and let $\left\lbrace {\mathcal{SOP}}(F_n, C, D_n) \right\rbrace$ be a sequence of approximated set-valued optimization problems. Assume that
\begin{enumerate}
\item[$(a)$] $F_n \underset{D_n, D}{\overset{\Gamma^{C}_{\text{seq}}}{\longrightarrow}} F$ $;$
\item[$(b)$] $\Ge(F_n, C, D_n) \neq \emptyset$ for all $n \in \N$ and Hypothesis $(H)$ is satisfied for $\mathcal{T} = \mathcal{G}$.
\end{enumerate}
Then, for all $\bx \in \Ge(F,C,D)$, there is a sequence $\dis \left\lbrace x_n \right\rbrace \subset X$ with $x_n \in   \Ge(F_n, C, D_n)$ for all $n \in \N$ and a pair $(n_k, \bz) \in \Nss \times \Ge(F, C, D)$ such that $x_{n_{k}} \stackrel{\Vert \cdot \Vert_{X}}{\longrightarrow} \bz$ and $F(\bz) \stackrel{l}{\approx} F(\bx).$
\end{theorem}

\begin{proof} Let $\bx \in \Ge(F,C,D).$ We fix $\varepsilon \in \inte(C)$. Since $F_n \underset{D_n, D}{\overset{\Gamma^{C}_{\text{seq}}}{\longrightarrow}} F,$ there exists a sequence $\left\lbrace \varphi_n \right\rbrace \subset X$ with $\varphi \in D_n$, for all $n \in \N$ such that $\varphi_n \stackrel{\Vert \cdot \Vert_{X}}{\longrightarrow} \bx$ and $F_n(\varphi_n) - \varepsilon \; \lell F(\bx),$ eventually. According to $(b)$, we can construct a sequence $\left\lbrace x_n \right\rbrace \subset X$ with $x_n \in \Ge(F_n,C,D_n) \cap  {\rm Lev}_{F_{n}(\varphi_{n})}(F_n)$ for all $n \in \N.$ Thus, we obtain that 
\begin{equation}\label{EqLowerGeofStab1}
F_n(x_n) - \frac{\varepsilon}{4} \; \lell F(\bx), \; \text{eventually.}
\end{equation}
Let us observe that $x_n \in D_n$ for all $n \in \N$. Because  $F_n \underset{D_n, D}{\overset{\Gamma^{C}_{\text{seq}}}{\longrightarrow}} F,$ then $\left( \left\lbrace D_n \right\rbrace, D \right)$ is a Kuratowski pair. Hence, we can find a pair $(n_k, \bz) \in \Nss \times D$ such that $x_{n_{k}} \stackrel{\Vert \cdot \Vert_{X}}{\longrightarrow} \bz.$ Invoking $(a)$ again, we deduce that 
\begin{equation}\label{EqLowerGeofStab2}
F(\bz) - \frac{\varepsilon}{4} \; \lell F_{n_{k}}(x_{n_{k}}), \; \text{ eventually.}
\end{equation}
Combining \eqref{EqLowerGeofStab1} and \eqref{EqLowerGeofStab2}, we obtain from Lemma \ref{lemmastrictStrong} $(a)$ that $F(\bz) - \varepsilon \; \ll F(\bx).$ As $\varepsilon$ has been chosen arbitrarily, we deduce from Lemma \ref{lemmastrictStrong} $(b)$ that $F(\bz) \; \lell F(\bx).$ Finally, from the very nature of $\bx$, we get $F(\bz) \stackrel{l}{\approx} F(\bx).$
\end{proof}

\noindent Finally, we consider the case of relaxed minimal solutions. We obtained the following result. 

\begin{theorem}[Internal stability of relaxed minimal solutions]\label{thLowerStabRelax} Consider the set-valued optimization problem ${\mathcal{SOP}}(F, C, D)$ and let $\left\lbrace {\mathcal{SOP}}(F_n, C, D_n) \right\rbrace$ be a sequence of approximated set-valued optimization problems. Assume that
\begin{enumerate}
\item[$(a)$] $F_n \underset{D_n, D}{\overset{\Gamma^{C}_{\text{seq}}}{\longrightarrow}} F$ $;$
\item[$(b)$] $\Re(F_n, C, D_n) \neq \emptyset$ for all $n \in \N$ and Hypothesis $(H)$ is satisfied for $\mathcal{T} = \mathcal{R}$.
\end{enumerate}
Then, for all $\bx \in \Re(F,C,D)$, there is a sequence $\dis \left\lbrace x_n \right\rbrace \subset X$ with $x_n \in   \Re(F_n, C, D_n)$ for all $n \in \N$ and a pair $(n_k, \bz) \in \Nss \times \Re(F, C, D)$ such that $x_{n_{k}} \stackrel{\Vert \cdot \Vert_{X}}{\longrightarrow} \bz$ and $F(\bz) \; \lell F(\bx).$
\end{theorem}

\begin{proof} Let $\bx \in \Re(F,C,D).$ By adapting the method developed to prove Theorem \ref{TheoLowGeofStab}, we show that there is a sequence $\dis \left\lbrace x_n \right\rbrace \subset X$ with $x_n \in   \Re(F_n, C, D_n)$ for all $n \in \N$ and a pair $(n_k, \bz) \in \Nss \times D$ such that $x_{n_{k}} \stackrel{\Vert \cdot \Vert_{X}}{\longrightarrow} \bz$ and $F(\bz) \; \lell F(\bx).$ It remains to show that $\bz \in \Re(F,C,D).$ Assume to the contrary that $\bz \notin \Re(F,C,D).$ Then, according to Remark \ref{RemarkLinkRelations}, we can find $x_0 \in D$ such that $F(x_0) \; \ll F(\bx).$ However, this yields to a contradiction as $F \in \F$ (see Proposition \ref{propcproper}).
\end{proof}

\begin{remark} From Remark \ref{RemInc}, Theorems \ref{TheoUppGeofStab} and \ref{TheoLowGeofStab} are valid for 
$\Se(F,C,D)$ and $\Pe(F,C,D)$ when $F$ is $C$-closed valued. 
\end{remark}

\noindent For the purposes of completeness, we illustrate the practical aspect of the tools introduced in our work. 

\begin{example} Let $X = Z = \R$ and let $C = \R_+$. We set $D =  \left[0, \dis \frac{\pi}{4}\right]$ and $D_n =  \left[0, \dis \frac{\pi}{4} + \frac{\pi}{n} \right),$ for all $n \in \N.$ Consider the following set-valued optimization problem:
\begin{equation*}
{\mathcal{SOP}}\left(F, \R_+, \left[0, \dis \frac{\pi}{4}\right] \right) : \left\lbrace
\begin{array}{lc}
 {\rm Minimize~} & F(x) := \left( \sin(x) , 3 \right),\\ 
 {\rm with~respect~to} & \R_+,\\ 
{\rm and~subject~to~}  & x \in \left[0, \dis \frac{\pi}{4}\right],
\end{array}
\right.
\end{equation*}
along with the sequence of perturbed set-valued optimization problems defined by:
\begin{equation*}
{\mathcal{SOP}}\left(F_n, \R_+, \left[0, \dis \frac{\pi}{4} + \frac{\pi}{n} \right) \right) : \left\lbrace
\begin{array}{lc}
 {\rm Minimize~} & F_n (x) := \left( \sin(x\left(1 + \dis \frac{1}{n+1}\right)) , \; 3 + \exp (n) \right),\\ 
 {\rm with~respect~to} & \R_+,\\ 
{\rm and~subject~to~}  & x \in \left[0, \dis \frac{\pi}{4} + \frac{\pi}{n + 1} \right), \; n \in \N.
\end{array}
\right.
\end{equation*}
First, observe that $\Re(F,C,D) = \left\lbrace 0 \right\rbrace$. We show the external and the internal stability of the set of relaxed minimal solutions to ${\mathcal{SOP}}\left(F, \R_+, \left[0, \dis \frac{\pi}{4}\right]\right).$ We proceed in two steps. \\
\textbf{Step 1.} We show that $\left( \left\lbrace D_n \right\rbrace, D \right)$ is a Kuratowski pair. Note that $D$ is  closed. Invoking the hit-and-miss criterion $($see Rockafellar and Wets \cite[\text{Theorem 4.5}]{RWE98} $)$ we deduce that $D_n \stackrel{PK}{\rightharpoonup} D$. Moreover, for all $n \geq 1$, $D_n \subset D_1$ and ${\rm cl}(D_1)$ is compact. Thus, it follows from Proposition \ref{criterionKuratowskiPair} that $\left( \left\lbrace D_n \right\rbrace, D \right)$ is a Kuratowski pair.\\
\textbf{Step 2.} We prove that $F_n \underset{D_n, D}{\overset{\Gamma^{C}_{\text{seq}}}{\longrightarrow}} F$. Let $\varepsilon \in \inte(C) = \left(0, + \infty \right)$ and $\bx \in D.$ From the continuity of $\sin(\cdot)$ on $\R$, 
we deduce that for any sequence $\left\lbrace x_n \right\rbrace \subset \R$ with $x_n \in D_n$ for all $n \in \N$ such that $x_{n} \stackrel{\vert \cdot \vert}{\longrightarrow} \bx,$ one has  $\dis \sin\left(x_n \left(1 + \frac{1}{n + 1}\right) \right) \in \left( \sin(\bx) - \varepsilon, \sin(\bx) + \varepsilon \right),$ eventually. This yields that 
\begin{equation*}
F_n(x_n) := \left[\sin\left(x_n \left(1 + \frac{1}{n+1}\right) \right) , 3 + \exp(n) \right) \subset \left[ \sin(\bx) - \varepsilon, + \infty\right) =   {\rm cl}\left(F(\bx) - \varepsilon + C\right),
\end{equation*} 
for $n$ large enough. In other words, $F(\bx) - \varepsilon \; \lell F_n(x_n),$ eventually. On the other hand, we consider the sequence $\left\lbrace x_n^\star \right\rbrace \subset \R$ define for all $n \in \N$ by $x_n^\star = 0.$ Thus, $x_n^\star \in D_n$ for all $n \in \N$ and $\sin(x_n^\star) \leq \sin(\bx)$. Hence, we get 
\begin{equation*}
F(\bx) + \varepsilon := \left( \sin(\bx) + \varepsilon, 3 + \varepsilon \right) \subset \left[ \sin(0), + \infty\right) = {\rm cl}\left(F_n(x^\star_n) + C\right), \; \text{ for all } n \in \N.
\end{equation*} 
In other words, $F_n(x_n^\star) \; \lell F(\bx) + \varepsilon,$ eventually. Finally,  $F_n \underset{D_n, D}{\overset{\Gamma^{C}_{\text{seq}}}{\longrightarrow}} F.$ Thus, we can apply Theorem \ref{UpperStabilityRelaxed}. Moreover, since Hypothesis $(H)$ is satisfied for $\mathcal{T} = \mathcal{R}$, we can also apply Theorem \ref{thLowerStabRelax}.
\end{example}

\section{Conclusion and open questions}\label{Section 4}

Two variational convergences, namely the Gamma-cone convergence and the sequential Gamma-cone convergence, are introduced in this work in order to study different kind of stability problems in Geoffroy spaces. In particular, we obtain several results related to the asymptotic properties of strong level sets and to the external and internal stability of minimal solutions to set-valued optimization problems. Our work enable to consider that both the objective function and the admissible domain are perturbed.\\

In Definition \ref{DefGeoffroySol}, we have introduced the concept of Geoffroy minimal solutions which are weaker than Pareto minimal solutions. The existence of Pareto minimal solutions to ${\mathcal{SOP}}(F, C, D)$ has been deeply studied by several authors (see, for instance, \cite{ARO05, HRM07, ZHK24}). In particular, Hern\'andez and Rodr\'iguez-M\'arin \cite[Theorem 5.9]{HRM07} proved the next result. 

\begin{theorem}[Sufficient condition for the existence of Pareto minimal solutions]\label{ExistenceHerRoMa} Assume that $D$ is compact. If $\mathcal{L}(y) := \left\lbrace x \in  D \mid F(x) \; \lel \left\lbrace y \right\rbrace \right\rbrace$ is closed for all $y \in F(D)$, then $\Pe(F,C,D)$ is nonempty.
\end{theorem}

\noindent Therefore, we are able to prove an existence result for Pareto minimal solutions to ${\mathcal{SOP}}(F, C, D)$ in Geoffroy space as shown below.

\begin{theorem}[Existence of Pareto minimal solutions]\label{TheoExistenceParetoLower} Let $F \in \FF$ and let $D \subset X$ be a compact set. If $F$ is lower $\topo$-semicontinuous and $C$-closed-valued on $D$, then $\Pe(F,C,D)$ is nonempty.
\end{theorem}

\begin{proof} Consider $y \in F(D)$ and let $\left\lbrace x_n \right\rbrace \subset \mathcal{L}(y)$ be a sequence such that $x_n \stackrel{\Vert \cdot \Vert_{X}}{\longrightarrow} \bx \in X$. According to Theorem \ref{ExistenceHerRoMa}, it is sufficient to prove that $\bx \in \mathcal{L}(y).$ Since $\mathcal{L}(y) \subseteq D$, we deduce that $\bx \in D.$ Hence, from the lower $\topo$-semicontinuity of $F$ on $D$ and Proposition \ref{propcharlsc}, we deduce that $F$ is sequentially lower $\topo$-semicontinuous at $\bx$, i.e. for all $\varepsilon \in \inte(C)$
\begin{equation}\label{EqExistence1}
F(\bx) - \varepsilon \; \ll F(x_n), \; \text{ for } n \text{ large enough}.
\end{equation}
Combining estimate \eqref{EqExistence1} with the fact that $\left\lbrace x_n \right\rbrace \subset \mathcal{L}(y)$, we obtain 
\begin{equation}\label{EqExistence2}
\text{ for all } \varepsilon \in \inte(C), \;  F(\bx) - \varepsilon \; \ll \left\lbrace y \right\rbrace.
\end{equation}
Thus, because $F(\bx)$ is $C$-closed, we deduce from Lemma \ref{lemmastrictStrong} $(b)$ that $\bx \in  \mathcal{L}(y)$ and $\Pe(F,C,D) \neq \emptyset$.
\end{proof}

\noindent Thus, it would be interesting to see if the Theory developped by Hern\'andez and Rodr\'iguez-M\'arin \cite{HRM07} can be adapted to our new preorder $\lell.$ That way, one can discuss the nonemptyness of $\Ge(F,C,D).$\\

\noindent On the other hand, it can be seen that Hypothesis $(H)$ is a key concept in our work. However, it is not certain that this hypothesis is always true. Shedding light on this question would be a major step forward.\\

\noindent Finally, discussing the case in which the cone $C$ is also pertubed would be interesting. More precisely, the external and internal stability of set-valued optimization problems should be explored via sequences of approximate set-valued optimization problems of the form
\begin{equation*}
{\mathcal{SOP}}(F_n, C_n, D_n) : \left\lbrace
\begin{array}{lc}
 {\rm Minimize~} & F_n(x),\\ 
 {\rm with~respect~to} & C_n,\\ 
{\rm and~subject~to~}  & x\in D_n.
\end{array}
\right.
\end{equation*}
\noindent However, it should be pointed out that the set order topology $\topo$ is intrinsically connected to the cone $C$. In fact, we prove the following result:

\begin{theorem}\label{TheoFiness} Let $C_1, C_2 \subsetneq Z$ be two nonempty closed, convex and solid cones such that $C_1 \subseteq C_2$. Then, $\dis \tau^{C_1}$ is finer than $\dis \tau^{C_2}.$
\end{theorem}

\noindent To this end, we first state an intermediate result.

\begin{lemma}\label{LemmaVectorMino} Let $c \in C$ and $u \in \inte(C)$. Then, there is $N \in \mathbb{N}$ such that $\dis \frac{u}{2} - \frac{c}{N} \in C$. In particular, we have $\left\lbrace c \right\rbrace \ll \left\lbrace N u \right\rbrace$. 
\end{lemma}

\begin{proof} We set $c_n := \dis \frac{c}{n}$, for all $n \in \N$. Since $u \in \inte(C),$ we can find a radius $r > 0$ such that $\ball_r \left( \dis \frac{u}{2}\right) \subset C$. Furthermore, for $n$ large enough, we have $\dis \bigg\Vert \dis \frac{u}{2} - c_n - \frac{u}{2} \bigg\Vert = \dis \| c_n \| < r.$ As a result, we can find $N \in \N$ such that $\dis \frac{u}{2} - c_N \in C,$ which yields that
\begin{center}
$ \dis \left\lbrace 0_Z \right\rbrace \lel \dis \left\lbrace \frac{u}{2} - c_N \right\rbrace \; \Longleftrightarrow \;  \left\lbrace c_N \right\rbrace \lel \left\lbrace \frac{u}{2} \right\rbrace \; \Longrightarrow \;  \left\lbrace c \right\rbrace \ll \left\lbrace Nu \right\rbrace.$
\end{center}
The proof is complete.
\end{proof}

\noindent Then, we can demonstrate Theorem \ref{TheoFiness}.

\begin{proof}[Proof of Theorem \ref{TheoFiness}.] Let $\hX \in \pzs$. From \cite[Proposition 4.4]{GLA22}, the families 
\begin{center}
$\Vo_{C_{1}} (\hX) := \left\lbrace \left(\hX-u; \hX+u \right)_{C_{1}} \mid  u \in \inte(C_{1}) \right\rbrace$ and $\Vo_{C_{2}} (\hX) := \left\lbrace \left(\hX-u; \hX+u \right)_{C_{2}} \mid  u \in \inte(C_{2}) \right\rbrace$ 
\end{center}
are two bases of neighborhoods of $\hX$ for $\dis \tau^{C_1}$ and $\dis \tau^{C_2},$ respectively. We show that each elements of $\Vo_{C_{2}} (\hX)$ contains an element of $\Vo_{C_{1}} (\hX)$. Let $\overline{u}_{2} \in \inte(C_2)$. We consider the element $\left(\hX - \overline{u}_{2} ; \hX + \overline{u}_{2} \right)_{C_{2}}$ of $\Vo_{C_{2}} (\hX)$. 
Now, let us take $\overline{u}_{1} \in \inte(C_1)$ and let us set $\dis \overline{u}_{1}^n := \frac{\overline{u}_{1}}{n}$, for some $n \in \N$. Since $C_1 \subseteq C_2$, it is clear that for all $n \in \N$, $\overline{u}_{1}^n \in C_2$. Moreover, according to Lemma \ref{LemmaVectorMino}, there exists $N \in \N$ such that $\left\lbrace \overline{u}_{1}^N \right\rbrace \; \prec^{l}_{C_{2}} \; \left\lbrace \overline{u}_{2} \right\rbrace.$
Therefore, since $\hX \lel \hX$, we get 
\begin{center}
$\hX - \overline{u}_{2} \; \prec^{l}_{C_{2}} \hX - \overline{u}_{1}^N \;$ and $ \; \hX +  \overline{u}_{1}^N \; \prec^{l}_{C_{2}} \hX + \overline{u}_{2}.$
\end{center}
Since $\overline{u}_{1}^N \in \inte(C_1)$, we deduce that $\left(\hX - \overline{u}_{1}^N ; \hX + \overline{u}_{1}^N \right)_{C_{1}} \subseteq \left(\hX - \overline{u}_{2} ; \hX + \overline{u}_{2} \right)_{C_{2}}.$
\end{proof}

\noindent Thus, the study of perturbation of this type could lead to some extremely promising results, which would probably go beyond the scope of stability analysis. \\

\noindent \textbf{Disclosure of interest.} The author report there are no competing interests to declare.



\begin{thebibliography}{cc}
\bibitem{ARO05}{\sc Alonso, M., \& Rodr\'iguez-Mar\'in, L.} (2005). Set-relations and optimality conditions in set-valued maps. Nonlinear Analysis: Theory, Methods \& Applications, 63(8), 1167--1179.

\bibitem{ADH20}{\sc Anh, L. Q., Duy, T. Q., \& Hien, D. V.} (2020). Stability of efficient solutions to set optimization problems. Journal of Global Optimization, 78, 563-580.

\bibitem{ADH20b}{\sc Anh, L. Q., Duy, T. Q., Hien, D. V., Kuroiwa, D., \& Petrot, N.} (2020). Convergence of solutions to set optimization problems with the set less order relation. Journal of Optimization Theory and Applications, 185, 416-432.

\bibitem{AHS22}{\sc Ansari, Q.H., Hussain, N. \& Sharma, P.K.} (2022). Convergence of the solution sets for set optimization problems. J. Nonlinear Var. Anal., 6(3), 165-183.

\bibitem{BRA06}{\sc Braides, A.} (2006). A handbook of $\Gamma$-convergence. In Handbook of Differential Equations: stationary partial differential equations (Vol. 3, pp. 101-213). North-Holland.

\bibitem{GGJ17}{\sc Gaydu, M., Geoffroy, M. H., Jean-Alexis, C., \& Nedelcheva, D.} (2017). Stability of minimizers of set optimization problems. Positivity, 21, 127-141.


\bibitem{GMN17}{\sc Geoffroy, M. H., Marcelin, Y., \& Nedelcheva, D.} (2017). Convergence of relaxed minimizers in set optimization. Optimization Letters, 11(8), 1677-1690.

\bibitem{GEO19}{\sc Geoffroy, M. H.} (2019). A topological convergence on power sets well-suited for set optimization. Journal of Global Optimization, 73, 567--581.

\bibitem{GLA22}{\sc Geoffroy, M. H., \& Larrouy, J.} (2022). A New Topological Framework and Its Application to Well-Posedness in Set-Valued Optimization. Numerical Functional Analysis and Optimization, 43(16), 1848--1883.

\bibitem{GMM16}{\sc Guti\'errez, C., Miglierina, E., Molho, E., \& Novo, V.} (2016). Convergence of solutions of a set optimization problem in the image space. Journal of Optimization Theory and Applications, 170, 358-371.

\bibitem{HAM05}{\sc Hamel, A. H.} (2005). Variational principles on metric and uniform spaces. Habilitation thesis, University of Halle-Wittenberg, Germany.

\bibitem{HHE10}{\sc Hamel, A. H., \& Heyde, F.} (2010). Duality for set-valued measures of risk. SIAM J. Finan. Math. 1(1):66?95.

\bibitem{HHR11}{\sc Hamel, A. H., Heyde, F., \& Rudloff, B.} (2011). Set-valued risk measures for conical market models. Math. Financ. Econ. 5, no. 1, 1--28.

\bibitem{HAN22}{\sc Han, Y.} (2022). Painlev\'e-Kuratowski convergences of the solution sets for set optimization problems with cone-quasiconnectedness. Optimization, 71(7), 2185-2208.

\bibitem{HRM07}{\sc Hern\' andez, E., \& Rodr\' \i guez-Mar\' \i n, L.} (2007). Existence theorems for set optimization problems,  Nonlinear Anal. 67, no. 6, 1726--1736. Springer Optim. Appl., 35, Springer, New York, 2010.

\bibitem{HLO22}{\sc Hern\'andez, E., \& L\'opez, R.} (2022). On epi-convergence for set-valued maps. Optimization, 71(2), 403-417.

\bibitem{HLO22b}{\sc Hern\'andez, E., \& L\'opez, R.} (2022). Stability in Set-Valued Optimization Problems Using Asymptotic Analysis and Epi-Convergence. Applied Mathematics \& Optimization, 86(1), 7.

\bibitem{KLA19}{\sc Karuna, \& Lalitha, C. S.} (2019). External and internal stability in set optimization using gamma convergence. Carpathian Journal of Mathematics, 35(3), 393-406.

\bibitem{KLA22}{\sc Karuna, \& Lalitha, C. S.} (2022). Essential stability in set optimization. Optimization, 71(7), 1803-1816.

\bibitem{KLA24}{\sc Karuna, \& Lalitha, C. S.} (2022). Convergence, Scalarization and Continuity of Minimal Solutions in Set Optimization. Journal of the Operations Research Society of China, 12(3), 773-793.

\bibitem{KTZ15} {\sc Khan, A. A., Tammer, C., \& Z\u{a}linescu, C.} (2015). Set-valued optimization. An introduction with applications. Vector Optimization. Springer, Heidelberg. xxii+765 pp. 

\bibitem{KHO18}{\sc Khoshkhabar-amiranloo, S.} (2018). Stability of minimal solutions to parametric set optimization problems. Applicable Analysis, 97(14), 2510-2522. 

\bibitem{KUR97}{\sc Kuroiwa, D.}  (1997). Some Criteria in Set-valued Optimization. RIMS Kokyuroku, (985) : 171--176.

\bibitem{KUR98}{\sc Kuroiwa, D.} (1998). The natural criteria in set-valued optimization. Research on nonlinear analysis and convex analysis No. 1031, 85--90.

\bibitem{KUR01}{\sc Kuroiwa, D.} (2001). On set-valued optimization. Nonlinear Analysis 47, no. 2, 1395--1400.

\bibitem{ORO08}{\sc Oppezzi, P., \& Rossi, A.} (2008). A convergence for vector valued functions. Optimization, 57(3), 435-448.

\bibitem{ORO08b}{\sc Oppezzi, P., \& Rossi, A. M.} (2008). A convergence for infinite dimensional vector valued functions. Journal of Global Optimization, 42, 577-586.

\bibitem{RWE98}{\sc Rockafellar, R. T., \& Wets, R. J. B.} (2009). Variational analysis (Vol. 317). Springer Science \& Business Media.

\bibitem{TKU93}{\sc Tanaka, T., \& Kuroiwa, D.} (1993). The convexity of $A$ and $B$ assures $\inte(A) + B = \inte(A + B)$. Applied mathematics letters, 6(1), 83-86.

\bibitem{ZHK24}{\sc Zhou, Z., Huang, M., \& Kobis, E.} (2024). Globally proper efficiency of set optimization problems based on the certainly set less order relation. Applicable Analysis, 103(1), 184--197.

\end{thebibliography}
\end{document}